\documentclass[12pt,letter]{amsart}
\usepackage{amsmath,amssymb}
\title[Factors
with at most one Cartan subalgebra II]{On a class
of $\mathrm{II}_1$ factors\\
with at most one Cartan subalgebra II}
\author[N. Ozawa]{Narutaka OZAWA$^*$}
\address{Department of Mathematical Sciences,
University of Tokyo, Komaba, \mbox{153-8914}\\
\indent Department of Mathematics, UCLA, Los Angeles, CA90095}
\email{narutaka@ms.u-tokyo.ac.jp}
\thanks{${}^*$ Supported by NSF-Grant and Sloan Foundation}
\author[S. Popa]{Sorin Popa$^{**}$}
\address{Department of Mathematics,
UCLA, Los Angeles, CA90095}
\email{popa@math.ucla.edu}
\thanks{${}^{**}$ Supported by NSF-Grant}

\date{July 14, 2008}
\newtheorem{thm}{Theorem}[section]
\newtheorem{lem}[thm]{Lemma}
\newtheorem{prop}[thm]{Proposition}
\newtheorem{cor}[thm]{Corollary}
\newtheorem{thmA}{Theorem}

\newtheorem{cora}{Corollary}

\theoremstyle{definition}
\newtheorem*{defn}{Definition}
\newtheorem{rem}[thm]{Remark}
\numberwithin{equation}{section}
\newcommand{\R}{{\mathbb R}}
\newcommand{\C}{{\mathbb C}}
\newcommand{\N}{{\mathbb N}}
\newcommand{\Z}{{\mathbb Z}}
\newcommand{\Q}{{\mathbb Q}}
\newcommand{\M}{{\mathbb M}}
\newcommand{\F}{{\mathbb F}}

\newcommand{\B}{{\mathbb B}}
\newcommand{\D}{{\mathcal D}}
\newcommand{\GG}{{\mathcal G}}

\newcommand{\NN}{{\mathcal N}}
\newcommand{\Nor}{{\mathcal N}}
\newcommand{\Zt}{{\mathcal Z}}
\newcommand{\U}{{\mathcal U}}
\newcommand{\G}{\Gamma}
\newcommand{\e}{\varepsilon}
\newcommand{\op}{\mathrm{op}}
\newcommand{\p}{\varphi}
\newcommand{\hh}{{\mathcal H}}
\newcommand{\hk}{{\mathcal K}}

\DeclareMathOperator{\Ad}{Ad}
\DeclareMathOperator{\Aut}{Aut}
\DeclareMathOperator{\lspan}{span}
\newcommand{\vt}{\mathbin{\bar{\otimes}}}

\DeclareMathOperator*{\Lim}{Lim}
\newcommand{\id}{\mathrm{id}}
\DeclareMathOperator{\Tr}{Tr}
\DeclareMathOperator{\Ind}{Ind}
\DeclareMathOperator{\rank}{rank}
\newcommand{\SL}{\mathrm{SL}}
\newcommand{\PSL}{\mathrm{PSL}}
\newcommand{\SO}{\mathrm{SO}}
\newcommand{\SU}{\mathrm{SU}}
\newcommand{\ip}[1]{\mathopen{\langle}#1\mathclose{\rangle}}

\DeclareMathOperator{\dom}{dom}
\DeclareMathOperator{\ran}{ran}
\newcommand{\same}{\ \rule[3pt]{50pt}{0.5pt}\ }
\newcommand{\HH}{$(\mathrm{HH})$}
\begin{document}
\begin{abstract}
This is a continuation of our previous paper studying the structure
of Cartan subalgebras of von Neumann factors of type $\mathrm{II}_1$.
We provide more examples of $\mathrm{II}_1$ factors having either zero,
one or several Cartan subalgebras. We also prove a rigidity result
for some group measure space $\mathrm{II}_1$ factors.
\end{abstract}
\maketitle
\section{Introduction}
A celebrated theorem of Connes (\cite{connes:cls}) states that all
amenable $\mathrm{II}_1$ factors are isomorphic to the approximately
finite dimensional $\mathrm{II}_1$ factor $R$ of Murray and von
Neumann. In particular, all group $\mathrm{II}_1$ factors
$L(\Gamma)$ associated with ICC (infinite conjugacy class) amenable
groups $\Gamma$, and all group measure space $\mathrm{II}_1$ factors
$L^\infty(X)\rtimes\Gamma$ arising from (essentially) free ergodic
probability-measure-preserving (abbreviated as p.m.p.) actions
$\Gamma \curvearrowright X$ of countable amenable groups $\Gamma$ on
standard probability spaces $X$, are isomorphic to $R$.

In contrast to the amenable case, the group measure space
$\mathrm{II}_1$ factors $L^\infty(X)\rtimes\G$ of free ergodic
p.m.p.\ actions of non-amenable groups $\G$ on standard probability
spaces $X$ form a rich and particularly important class of
$\mathrm{II}_1$ factors. More general crossed product construction
provides a wider class. We want to investigate the isomorphy problem
of the crossed product $\mathrm{II}_1$ factors. Namely, given the
crossed product $M=Q\rtimes\G$ of a finite amenable von Neumann
algebra $(Q,\tau)$ by a $\tau$-preserving action of a countable
group $\G$, to what extent can we recover information on the
original action $\G\curvearrowright Q$? In particular, does there
exist a group measure space $\mathrm{II}_1$ factor
$M=L^\infty(X)\rtimes\G$ which remembers completely the group $\G$
and the action $\G\curvearrowright X$? The first task would be to
determine all regular amenable subalgebras of a given
$\mathrm{II}_1$ factor $M$. Recall that a von Neumann subalgebra $P$
of $M$ is said to be \emph{regular} if the normalizer group of $P$
in $M$ generates $M$ as a von Neumann algebra (\cite{dixmier}). A
regular maximal abelian subalgebra $A$ of $M$ is called a
\emph{Cartan subalgebra} (\cite{feldman-moore}). In the case of a
group measure space $\mathrm{II}_1$ factor $M=L^\infty(X)\rtimes\G$,
the von Neumann subalgebra $L^\infty(X)$ is a Cartan subalgebra and
determining its position amounts to recovering the orbit equivalence
relation of the original action $\G\curvearrowright X$ (see
\cite{feldman-moore}). By \cite{cfw}, the approximately finite
dimensional $\mathrm{II}_1$ factor $R$ has a unique Cartan
subalgebra, up to conjugacy by an automorphism of $R$. In the
previous paper (\cite{I}), we provided the first class of examples
of non-amenable $\mathrm{II}_1$ factors having unique Cartan
subalgebra. They are the group measure space $\mathrm{II}_1$ factors
$M=L^\infty(X)\rtimes\F_r$ associated with free ergodic p.m.p.\
profinite actions $\F_r\curvearrowright X$ of free groups $\F_r$. In
this paper, we extend this result from the free groups $\F_r$ to a
lager class of countable groups with the property (strong)
{\HH}$^+$, defined as follows.

\begin{defn}
Let $G$ be a second countable locally compact group.
By a \emph{$1$-cocycle}, we mean a continuous map $b\colon G\to\hk$,
or a triplet $(b,\pi,\hk)$ of $b$ and a continuous unitary $G$-representation
$\pi$ on a Hilbert space $\hk$, which satisfies the $1$-cocycle identity:
\[
\forall g,h\in\G,\quad b(gh)=b(g)+\pi_gb(h).
\]
The $1$-cocycle $b$ is called \emph{proper} if the set $\{g\in G :
\|b(g)\|\le R\}$ is compact for every $R>0$. Assume that $G$ is
non-amenable. We say \emph{$G$ has the Haagerup property} (see
\cite{ccjjv,bo}) if it admits a proper $1$-cocycle $(b,\pi,\hk)$. In
the case when $\pi$ can be taken non-amenable (resp.\ to be weakly
contained in the regular representation), we say \emph{$G$ has the
property (resp.\ strong) {\HH}}. We say \emph{$G$ has the property
(strong) {\HH}$^+$} if $\G$ has the property (strong) {\HH} and the
complete metric approximation property (i.e., it is weakly amenable
with constant $1$).
\end{defn}

In Section \ref{sec:HH}, we will prove that
lattices of products of $\SO(n,1)$ $(n\geq2)$ and $\SU(n,1)$ have the property {\HH}$^+$,
and that lattices of $\SL(2,\R)$ and $\SL(2,\C)$ have the property strong {\HH}$^+$.
Building on our previous work (\cite{I}) and Peterson's deformation
technology (\cite{peterson}), we obtain the following.

\begin{thmA}\label{thm:A}
Let $M=Q\rtimes\G$ be the crossed product of a finite von Neumann algebra
$(Q,\tau)$ by a $\tau$-preserving action of a countable group $\G$
with the property {\HH}. Let $P\subset M$ be a regular weakly compact
von Neumann subalgebra. Then, $P\preceq_MQ$.
\end{thmA}

Since $L^\infty(X)\rtimes\G$ has the complete metric approximation
property if $\G$ has it and the action is profinite, the weak
compactness assumption holds automatically (see Section
\ref{sec:prelim}).

\begin{cora}\label{cor:A}
Let $\G$ be a countable group with the property {\HH}$^+$. Then,
$L(\G)$ has no Cartan subalgebra. Moreover, if
$\G\curvearrowright X$ is a free ergodic p.m.p.\ profinite action,
then $L^\infty(X)$ is the unique Cartan subalgebra in
$L^\infty(X)\rtimes\G$, up to unitary conjugacy.
\end{cora}

As in \cite{I}, a stronger result holds if $\G$ has the property
strong {\HH}.

\begin{thmA}\label{thm:B}
Let $M=Q\rtimes\G$ be the crossed product
of a finite amenable von Neumann algebra $(Q,\tau)$ by a $\tau$-preserving
action of a countable group $\G$ with the property strong {\HH}.
Let $P\subset M$ be an amenable von Neumann subalgebra such that
$P\not\preceq_MQ$ and $\GG\subset\Nor_M(P)$ be a subgroup
whose action on $P$ is weakly compact.
Then, the von Neumann subalgebra $\GG''$ is amenable.
\end{thmA}

\begin{cora}\label{cor:B}
Let $\G$ be a countable group with the property strong {\HH}$^+$.
Then, $L(\G)$ is strongly solid, i.e., the normalizer of every
amenable diffuse subalgebra generates an amenable von Neumann subalgebra.
\end{cora}

Once the Cartan subalgebra $L^\infty(X)$ is determined, the
isomorphy problem of $M=L^\infty(X)\rtimes\G$ reduces to that of the
orbit equivalence relations. Then, the group $\G$ and the action
$\G\curvearrowright X$ can be recovered if the orbit equivalence
cocycle untwists (\cite{zimmer}). Ioana (\cite{ioana:cocycle})
proved a cocycle (virtual) super-rigidity result with discrete
targets for p.m.p.\ profinite actions of property (T) groups. Here,
we prove a similar result for property $(\tau)$ groups, but with
some restrictions on the targets. Recall that a (residually finite)
group $\G$ is said to \emph{have the property $(\tau)$} if the
trivial representation is isolated among finite unitary
representations. See \cite{lubotzky,lubotzky-zuk} for more
information on this property.
\begin{thmA}\label{thm:D}
Let $\G=\G_1\times\G_2$ be a group with the property $(\tau)$
and $\G\curvearrowright X=\varprojlim X_n$ be a p.m.p.\ profinite action
with growth condition such that both $\G_i\curvearrowright X$
are ergodic. Let $\Lambda$ be a residually-finite group.
Then, any cocycle
\[
\alpha\colon\G\times X\to\Lambda
\]
virtually untwists, i.e., there exist $n\in\N$ and
a cocycle $\beta\colon\G\times X_n\to\Lambda$
which is equivalent to $\alpha$.
\end{thmA}

It is plausible that the residual finiteness assumption on $\Lambda$
is in fact redundant. Since there are groups having both properties
{\HH$^+$} and $(\tau)$, Theorems \ref{thm:A} and \ref{thm:D}
together imply a rigidity result for group measure space von Neumann
algebras. Let $\G'\le\G$ be a finite index subgroup and
$\G'\curvearrowright(X',\mu')$ be a m.p.\ action. Then, the induced
action $\Ind_{\G'}^{\G}(\G'\curvearrowright X')$ is the $\G$-action
on the measure space $\G/\G'\times X'$, given by
$g(p,x)=(gp,\sigma(gp)^{-1}g\sigma(p)x)$, where $\sigma$ is a fixed
cross section $\sigma\colon \G/\G'\to\G$. (The action is unique up
to conjugacy.) We say that two p.m.p.\ actions
$\G_i\curvearrowright(X_i,\mu_i)$, $i=1,2$, are \emph{strongly
virtually isomorphic} if there are a p.m.p.\ action
$\G'\curvearrowright(X',\mu')$ and finite index inclusions
$\G'\hookrightarrow\G_i$ such that $\G_i\curvearrowright X_i$ are
measure-preservingly conjugate to
$\Ind_{\G'}^{\G_i}(\G'\curvearrowright X')$.

\begin{cora}\label{cor:D}
Let $\G_i=\PSL(2,\Z[\sqrt{2}])$ and $p_1<p_2<\cdots$ be prime
numbers. Let $\G=\G_1\times\G_2$ act on
$X=\varprojlim\PSL(2,(\Z/p_1\cdots p_n\Z)[\sqrt{2}])$ by the
left-and-right translation action. Let $\Lambda\curvearrowright Y$
be any free ergodic p.m.p.\ action of a residually-finite group
$\Lambda$ and suppose that $L^\infty(X)\rtimes\G\cong
(L^\infty(Y)\rtimes\Lambda)^t$ for some $t>0$. Then, $t\in\Q$ and
the actions $\G\curvearrowright X$ and $\Lambda\curvearrowright Y$
are strongly virtually isomorphic.
\end{cora}

Connes and Jones (\cite{connes-jones}) gave a first example of
$\mathrm{II}_1$ factors having more than one Cartan subalgebra. We
present here a new class of examples. To describe it, recall first
that if $\G$ is a discrete group having an infinite normal abelian
subgroup $H$, then $L(H)$ is a Cartan subalgebra of $L(\G)$ if and
only if it satisfies the relative ICC condition: for any
$g\in\G\setminus H$, the set $\{ aga^{-1} : a\in H\}$ is infinite.
The group $H\rtimes\G$ acts on $H$ by $(a,g)b=agbg^{-1}$ (cf.\
Proposition 2.11 in \cite{I}).

\begin{thmA}\label{thm:C}
Let $\G\curvearrowright X$ be a free ergodic p.m.p.\ action of a
discrete group $\G$ having an infinite normal abelian subgroup $H$
satisfying the relative ICC condition. Assume that
$H\curvearrowright X$ is ergodic and profinite. Then, both
$L^\infty(X)$ and $L(H)$ are Cartan subalgebras of
$L^\infty(X)\rtimes\G$. Assume moreover that $\G\curvearrowright X$
is profinite and there is no $H\rtimes\G$-invariant mean on
$\ell^\infty(H)$. Then, the Cartan subalgebras $L^\infty(X)$ and
$L(H)$ are non-conjugate.
\end{thmA}

We distinguish two Cartan subalgebras by weak compactness.
The simplest example is the following. Another example will be
presented in Section \ref{sec:more}.

\begin{cora}\label{cor:c}
Let $p_1,p_2,\ldots$ be prime numbers.
Then the $\mathrm{II}_1$-factor
\[
M=L^\infty(\varprojlim(\Z/p_1\cdots p_n\Z)^2)\rtimes(\Z^2\rtimes\SL(2,\Z))
\]
has more than one Cartan subalgebra.
\end{cora}
We observe that in the above, $L(\Z^2)$ is actually an (strong) HT
Cartan subalgebra of $M$, in the sense of \cite{popa:betti}. Thus,
while an HT factor has unique HT Cartan subalgebra, up to unitary
conjugacy, there exist HT factors that have at least two
non-conjugate Cartan subalgebras. It is plausible that there is no
essentially-free group action which gives rise to the same orbit
equivalence relation as $(L(\Z^2)\subset M)$. Such examples were
first exhibited by Furman (\cite{furman}). See also \cite{ms} and
\cite{popa:cocycle}.

\section{Groups with the property {\HH}}\label{sec:HH}
Let $G$ be a locally compact group.
We recall that a unitary $\G$-representation $(\pi,\hh)$
is called \emph{amenable} if there is a state $\p$ on $\B(\hh)$
which is $\Ad\pi$-invariant: $\p\circ\Ad\pi_g=\p$ for all $g\in G$.
This notion was introduced and studied by Bekka (\cite{bekka}).
Among other things, he proved that $\pi$ is amenable if and only if
$\pi\otimes\bar{\pi}$ weakly contains the trivial representation.

Let $\sigma$ be the conjugate action of $G$ on $L^\infty(G)$:
$(\sigma_hf)(g)=f(h^{-1}gh)$ for $f\in L^\infty(G)$ and $g,h\in G$.
We say a locally compact group $G$ is \emph{inner-amenable} if
there is a $\sigma$-invariant state $\mu$ on $L^\infty(G)$ which
vanishes on $C_0(G)$. We note that in several literatures it is
only required that $\mu$ is $\sigma$-invariant and $\mu\neq\delta_e$
(in case $G$ is discrete).

\begin{prop}
A locally compact group $G$ with the property {\HH}
has the Haagerup property and is not inner-amenable.
\end{prop}
\begin{proof}
Let $(b,\pi,\hk)$ be a proper $1$-cocycle and suppose
that there is a singular $\sigma$-invariant state $\mu$
on $L^\infty(G)$.
For $x\in\B(\hh)$, we define $f_x\in L^\infty(G)$ by
$f_x(g)=\|b(g)\|^{-2}\ip{xb(g),b(g)}$.
Let $h\in G$ be fixed.
Since
\[
\|b(h^{-1}gh)-\pi_h^{-1}b(g)\|=\|b(h^{-1})+\pi_{h^{-1}g}b(h)\|\le2\|b(h)\|,
\]
and $\|b(g)\|\to\infty$ as $g\to\infty$,
one has $\sigma_h(f_x)-f_{\pi_hx\pi_h^*}\in C_0(G)$.
It follows that the state $\p$ on $\B(\hh)$ defined by
$\p(x)=\mu(f_x)$ is $\Ad\pi$-invariant.
This means $\pi$ is amenable.
\end{proof}

We do not know whether the converse is also true.
Combined with Proposition 2.11 in \cite{I}, the above proposition yields
the following.
\begin{cor}
A discrete group $\G$ with the property {\HH}$^+$
does not have an infinite normal amenable subgroup.
\end{cor}

We are indebted to Y. Shalom and Y. de Cornulier respectively for
$(\ref{item:shalom})$ and $(\ref{item:cornulier})$ of the next
statement.

\begin{thm}\label{thm:examplehh}
The following are true.
\begin{enumerate}
\item
Each of the properties {\HH}, {\HH}$^+$, strong {\HH} and strong {\HH}$^+$
inherits to a lattice of a locally compact group.
\item
If $G_1$ and $G_2$ have the property {\HH} (resp.\ {\HH}$^+$),
then so does $G_1\times G_2$.
\item\label{item:shalom}
The groups $\SO(n,1)$ with $n\geq2$ and $\SU(n,1)$
have the property {\HH}$^+$.
\item\label{item:cornulier}
The groups $\SL(2,\R)$ and $\SL(2,\C)$ have
the property strong {\HH}$^+$.
\item
Suppose $\G$ is a countable non-amenable group acting properly on
a finite-dimensional $\mathrm{CAT}(0)$ cube complex.
If all hyperplane stabilizer groups are non-co-amenable,
then $\G$ has the property {\HH}$^+$.
If all hyperplane stabilizer groups are amenable,
then $\G$ has the property strong {\HH}$^+$.
\end{enumerate}
\end{thm}
\begin{proof}
The assertion $(1)$ follows from the fact that
the restriction of non-amenable (resp.\ weakly sub-regular)
representation to a lattice is non-amenable (resp.\ weakly sub-regular).
For the assertion $(2)$, just consider the direct sum of $1$-cocycles.
We prove the property {\HH}$^+$ for $G=\SO(n,1)$ $(n\geq2)$ and $\SU(n,1)$.
It follows from Theorem 3 in \cite{bv} that every non-trivial irreducible
representation of $G$ is non-amenable.
Since $G$ does not have the property (T),
by \cite{shalom1}, there is a non-trivial irreducible representation
with an unbounded $1$-cocycle.
But, by \cite{shalom2}, every unbounded $1$-cocycles of
$G$ is proper.
Thus, $G$ has the property {\HH}.
Weak amenability is proved in \cite{decanniere-haagerup,cowling}.

The irreducible representation of
$\SL(2,\R)$ and $\SL(2,\C)$ which have non-trivial $1$-cocycles
are found in the principal series (see Example 3 in \cite{guichardet})
and hence are weakly equivalent
to the regular representation.

If a group $\G$ acts properly on a $\mathrm{CAT}(0)$ cube complex $\Sigma$,
then it has a proper $1$-cocycle into the $\ell^2(H)$, where $H$
is the set of hyperplanes in $\Sigma$. (See \cite{niblo-reeves}.)
The unitary representation on $\ell^2(H)$ is non-amenable (resp.\ weakly
contained in the regular representation) if and only if
all hyperplane stabilizer subgroups are non-co-amenable (resp.\ amenable).
Weak amenability for finite-dimensional $\mathrm{CAT}(0)$ cube complexes
is proved in \cite{guentner-higson,mizuta}.
\end{proof}

Note that by a result of \cite{csv}, the wreath product
$(\Z/2\Z)\wr\F_2$ acts properly on an infinite-dimensional
$\mathrm{CAT}(0)$ cube complex with all hyperplane stabilizer
subgroups amenable (being subgroups of $\bigoplus_{\F_2}\Z/2\Z$). It
follows that $(\Z/2\Z)\wr\F_2$ has the property strong {\HH}, but
not {\HH}$^+$.
\section{Miscellaneous Results}\label{sec:prelim}
We use the same conventions and notations as in the previous paper
(\cite{I}). Thus the symbol ``$\Lim$'' will be used for a state on
$\ell^\infty(\N)$, or more generally on $\ell^\infty(I)$ with $I$
directed, which extends the ordinary limit, and that the
abbreviation ``u.c.p.''\ stands for ``unital completely positive.''
We say a map is \emph{normal} if it is ultraweakly continuous.
Whenever a \emph{finite} von Neumann algebra $M$ is being
considered, it comes equipped with a distinguished faithful normal
tracial state, denoted by $\tau$. Any group action on a finite von
Neumann algebra is assumed to preserve the tracial state $\tau$. If
$M=P\rtimes\G$ is a crossed product von Neumann algebra, then the
tracial state $\tau$ on $M$ is given by
$\tau(au_g)=\delta_{g,e}\tau(a)$ for $a\in P$ and $g\in\G$. A von
Neumann subalgebra $P\subset M$ inherits the tracial state $\tau$
from $M$, and the unique $\tau$-preserving conditional expectation
from $M$ onto $P$ is denoted by $E_P$. We denote by $\Zt(M)$ the
center of $M$; by $\U(M)$ the group of unitary elements in $M$; and
by
\[
\Nor_M(P)=\{ u\in\U(M) : (\Ad u)(P)=P\}
\]
the normalizer group of $P$ in $M$, where $(\Ad u)(x)=uxu^*$. A von
Neumann subalgebra $P\subset M$ is called \emph{regular} if
$\Nor_M(P)''=M$. A regular maximal abelian von Neumann subalgebra
$A\subset M$ is called a \emph{Cartan subalgebra}. We note that if
$\G\curvearrowright X$ is a free ergodic p.m.p.\ action, then
$A=L^\infty(X)$ is a Cartan subalgebra in the crossed product
$L^\infty(X)\rtimes\G$. (See \cite{feldman-moore}.)

We recall the definition of weak compactness.
\begin{defn}
Let $(P,\tau)$ be a finite von Neumann algebra,
and $\G\curvearrowright P$ be a $\tau$-preserving action.
The action is called \emph{weakly compact} if
there is a net $\eta_n\in L^2(P\vt\bar{P})_+$ such that
\begin{enumerate}
\item $\|\eta_n-(v\otimes\bar{v})\eta_n\|_2\to0$ for $v\in\U(P)$;
\item $\|\eta_n-\Ad(u\otimes\bar{u})(\eta_n)\|_2\to0$ for $u\in\G$;
\item $\ip{(a\otimes1)\eta_n,\eta_n}=\tau(a)$ for all $a\in P$.
\end{enumerate}
(These conditions force $P$ to be amenable.)
A von Neumann subalgebra $P$ of $M$ is called \emph{weakly compact} if
the action $\Nor_M(P)\curvearrowright P$ is weakly compact.
\end{defn}

It is proved in \cite[Proposition 3.4]{I} that if
$\G\curvearrowright Q$ is weakly compact,
then $Q$ is weakly compact in the crossed product $Q\rtimes\G$.

\begin{thm}[Theorem 3.5 in \cite{I}]\label{thm:cmapwc}
Let $M$ be a finite von Neumann algebra with
the complete metric approximation property.
Then, every amenable von Neumann subalgebra
$P$ is weakly compact in $M$.
\end{thm}

Let $Q\subset M$ be finite von Neumann algebras.
Then, the conditional expectation $E_Q$ can be viewed
as the orthogonal projection $e_Q$ from $L^2(M)$ onto
$L^2(Q)\subset L^2(M)$.
It satisfies $e_Qxe_Q=E_Q(x)e_Q$ for every $x\in M$.
The \emph{basic construction} $\ip{M,e_Q}$ is the
von Neumann subalgebra of $\B(L^2(M))$ generated by
$M$ and $e_Q$. We note that $\ip{M,e_Q}$ coincides
with the commutant of the right $Q$-action in $\B(L^2(M))$.
The conditional expectation $E_Q$ extends on $\ip{M,e_Q}$
by the formula $E_Q(z)e_Q=e_Qze_Q$ for $z\in\ip{M,e_Q}$.
The basic construction
$\ip{M,e_Q}$ comes together with the faithful normal semi-finite
trace $\Tr$ such that $\Tr(xe_Qy)=\tau(xy)$.
We denote
\[
C^*(Me_QM)=\mbox{the norm-closed linear span of }\{ xe_Qy : x,y\in M\},
\]
which is an ultraweakly dense C$^*$-subalgebra of $\ip{M,e_Q}$.
Suppose that $\theta$ is a $\tau$-preserving u.c.p.\ map on $M$ such that
$\theta|_Q=\id_Q$. Then, $\theta$ can be regarded as a contraction
on $L^2(M)$ which commutes the left and right $Q$-actions on $L^2(M)$.
In particular, $\theta\in\ip{M,e_Q}$.
See Section 1.3 in \cite{popa:betti} for more information on
the basic construction.

The following is Theorem A.1 in \cite{popa:betti}
(see also Theorem 2.1 in \cite{popa:strong}).

\begin{thm}\label{finitebimodule}
Let $P,Q\subset M$ be finite von Neumann subalgebras.
Then, the following are equivalent.
\begin{enumerate}
\item
There exists a non-zero projection $e\in P'\cap\ip{M,e_Q}$
such that $\Tr(e)<\infty$.
\item
There exist non-zero projections
$p\in P$ and $q\in Q$, a normal $*$-homomorphism
$\theta\colon pPp\to qQq$ and a non-zero partial isometry $v\in M$
such that
\[
\forall x\in pPp\quad xv=v\theta(x)
\]
and $v^*v\in\theta(pPp)'\cap qMq$, $vv^*\in p(P'\cap M)p$.
\end{enumerate}
\end{thm}

\begin{defn}
Let $P,Q\subset M$ be finite von Neumann algebras. Following
\cite{popa:strong}, we say that \emph{$P$ embeds into $Q$ inside
$M$}, and write $P\preceq_MQ$, if any of the conditions in
Theorem~\ref{finitebimodule} holds.
\end{defn}

Let $P\subset\NN$ be von Neumann algebras.
We say a state $\p$ on $\NN$ is $P$-central if
$\p(u^*xu)=\p(x)$ for all $u\in\U(P)$ and $x\in\NN$,
or equivalently $\p(ax)=\p(xa)$ for all
$a\in P$ and $x\in\NN$.

\begin{lem}\label{lem:cptvsht}
Let $P,Q\subset M$ be a finite von Neumann algebras,
and $\p$ be a $P$-central state on $\ip{M,e_Q}$ whose restriction
to $M$ is normal.
If $P\not\preceq_M Q$, then $\p$ vanishes on $C^*(Me_QM)$.
\end{lem}
\begin{proof}
We assume $\p(C^*(Me_QM))\neq\{0\}$ and prove $P\preceq_M Q$.
Since $M$ sits inside the multiplier of $C^*(Me_QM)$,
there is an approximate unit $f_n$ of $C^*(Me_QM)$
such that $\|[f_n,u]\|\to0$ for every $u\in\U(M)$.
(That is, $(f_n)$ is a quasi-central approximate unit
for the ideal $C^*(Me_QM)$ in the C$^*$-algebra
$M+C^*(Me_QM)$.)
We may assume that each $f_n$ belongs to the linear span
of $\{xe_Qy : x,y\in M\}$.
We define positive linear functionals $\p_n$ and $\psi$
on $\ip{M,e_Q}$ by $\p_n(z)=\p(f_nzf_n)$ and
$\psi(z)=\Lim\p_n(z)$ for $z\in\ip{M,e_Q}$.
We note that $\psi$ is non-zero and still $P$-central.
We claim that $\psi$ is normal.
We observe that the net $(\p_n)$ actually norm converges to $\psi$
(w.r.t.\ $\Lim$), since $\Lim\psi(f_n)=\Lim\p(f_n)=\|\psi\|$.
Hence, it suffices to show that each $\p_n$ is normal.
Now let $n$ be fixed and $f_n=\sum_{i=1}^kx_ie_Qy_i$.
Then, for any $z\in\ip{M,e_Q}_+$, one has
\[
f_n^*zf_n=\sum_{i,j=1}^ky_i^*E_Q(x_i^*zx_j)e_Qy_j
\le\sum_{i,j=1}^ky_i^*E_Q(x_i^*zx_j)y_j\in M
\]
since $[E_Q(x_i^*zx_j)]_{i,j=1}^k$ is a positive element in $\M_k(Q)$
which commutes with $\mathrm{diag}(e_Q,\ldots,e_Q)$.
Hence, one has
\[
\p_n(z)=\p(f_n^*zf_n)\le(\p|_M)(\sum_{i,j=1}^ky_i^*E_Q(x_i^*zx_j)y_j).
\]
This implies that $\p_n$ is normal, and thus so is $\psi$.
It follows that $\psi$ can be regarded as a positive non-zero element
in $P'\cap L^1\ip{M,e_Q}$ (see Section IX.2 in \cite{takesaki:II}).
Taking a suitable spectral projection, we are done.
\end{proof}

We recall that A.1 in \cite{popa:betti} shows the following:
\begin{lem}\label{conjugatecartan}
Let $A$ and $B$ be Cartan subalgebras of
a type $\mathrm{II}_1$-factor $M$.
If $A\preceq_M B$, then there exists $u\in\U(M)$ such that $uAu^*=B$.
\end{lem}

Finally, we state some elementary lemmas about u.c.p.\ maps
and positive linear functionals.
We include proofs for the reader's convenience.
\begin{lem}\label{lem:appcs}
Let $(M,\tau)$ be a finite von Neumann algebra and
$\theta$ be a $\tau$-symmetric u.c.p.\ map on $M$.
Then for every $a,x\in M$, one has
\[
\|\theta(ax)-\theta(a)\theta(x)\|_2
\le2\|x\|_\infty\|a\|_\infty^{1/2}\|a-\theta(a)\|_2^{1/2}.
\]
\end{lem}
\begin{proof}
Let $\theta(x)=V^*\pi(x)V$ be a Stinespring dilation. Then,
\begin{align*}
\|\theta(ax)-\theta(a)\theta(x)\|_2
&=\|V^*\pi(x^*)(1-VV^*)\pi(a^*)V\widehat{1}\|_2\\
&\le\|x\|_\infty\|(1-VV^*)^{1/2}\pi(a^*)V\widehat{1}\|_2\\
&=\|x\|_\infty\tau\bigl(\theta(aa^*)-\theta(a)\theta(a^*)\bigr)^{1/2}.
\end{align*}
Since $\tau\circ\theta=\tau$, this completes the proof.
\end{proof}

\begin{lem}\label{lem:stateapprox}
Let $\p$ and $\psi$ be positive linear functional
on a C$^*$-algebra and $\e>0$.
Suppose that $\p(1)\geq\psi(1)$ and $\p(x)-\psi(x)\le\e\|x\|$
for all $x\geq0$. Then, one has $\|\p-\psi\|\le2\e$.
\end{lem}
\begin{proof}
Let $\p-\psi=(\p-\psi)_+-(\p-\psi)_-$ be the Hahn decomposition.
Since $(\p-\psi)(1)\geq0$, one has $\|(\p-\psi)_-\|\le\|(\p-\psi)_+\|\le\e$.
\end{proof}
\section{Peterson's Deformation}
We review the work of Peterson on
real closable derivations, in order to give a qualitative
version of Lemma 2.3 in \cite{peterson}.

Let $(M,\tau)$ be a finite von Neumann algebra.
An \emph{$M$-$M$ bimodule} is a Hilbert space $\hh$ together with
normal representations $\lambda$ of $M$ and $\rho$ of $M^{\op}$
such that $\lambda(M)\subset\rho(M^{\op})'$.
The action of $M$ is referred to as the left $M$-action and
the action of $M^{\op}$ is referred to as the right $M$-action.
We write intuitively $a\xi b$ for $\lambda(a)\rho(b^{\op})\xi$.
By a \emph{closable derivation}, we mean a map $\delta$ from
a weakly dense $*$-subalgebra $\D$ of $M$ into an $M$-$M$ bimodule $\hh$,
which is closable as an operator from $L^2(M)$ into $\hh$ and
satisfies the Leibniz's rule:
\[
\delta(xy)=\delta(x)y+x\delta(y)
\]
for every $x,y\in\D$.
Moreover, a derivation is always assumed to be \emph{real}:
there is a conjugate-linear isometric involution $J$ on $\hh$
such that $J(x\delta(y)z)=z^*\delta(y^*)x^*$ for every
$x,y,z\in\D$ (which is equivalent to another definition:
$\ip{\delta(x),\delta(y)z}=\ip{z^*\delta(y^*),\delta(x^*)}$
for every $x,y,z\in\D$).

Let $\hh$ be an $M$-$M$ bimodule and $\delta\colon M\to\hh$ be
a closable derivation whose closure is denoted by $\bar{\delta}$.
Thanks to the important work of \cite{davies-lindsay,sauvageot:quantum},
$\dom\bar{\delta}\cap M$ is still a weakly dense $*$-subalgebra
and $\bar{\delta}$ satisfies the Leibniz's rule there.
Hence, for notational simplicity, the closure $\bar{\delta}$
will be written as $\delta$.
We recycle some notations from \cite{peterson}:
\[
\Delta=\delta^*\delta,\
\zeta_\alpha=\sqrt{\frac{\alpha}{\alpha+\Delta}},\
\tilde{\delta}_\alpha=\alpha^{-1/2}\delta\circ\zeta_\alpha
\]
(note that $\ran\zeta_\alpha\subset\dom\Delta^{1/2}=\dom\delta$) and
\[
\tilde{\Delta}_\alpha=\alpha^{-1/2}\Delta^{1/2}\circ\zeta_\alpha
=\sqrt{\frac{\Delta}{\alpha+\Delta}},\
\theta_\alpha=1-\tilde{\Delta}_\alpha.
\]
All operators are firstly defined as Hilbert space operators.
Since $1-\sqrt{t}\le\sqrt{1-t}$ for all $0\le t\le1$, one has
$\theta_\alpha\le\zeta_\alpha$ and
\[
\|a-\zeta_\alpha(a)\|_2\le\|\tilde{\Delta}_\alpha(a)\|_2
=\|\tilde{\delta}_\alpha(a)\|_2\le\|a\|_2\le\|a\|_\infty
\]
for all $a\in M$.
By Lemma 2.2 in \cite{peterson},
the operators $\zeta_\alpha$ and $\theta_\alpha$
map $M\subset L^2(M)$ into $M$ and
are $\tau$-symmetric u.c.p.\ on $M$.

We recall from \cite{sauvageot:feller} the following facts:
$\psi_t=\exp(-t\Delta^{1/2})$ form a semigroup of u.c.p.\ maps on $M$.
Let
\begin{align*}
\Gamma(b^*,c)&=\Delta^{1/2}(b^*)c+b^*\Delta^{1/2}(c)-\Delta^{1/2}(b^*c)\\
&=\lim_{t\to 0}\frac{\psi_t(b^*c)-\psi_t(b^*)\psi_t(c)}{t}\in L^2(M)
\end{align*}
for $b,c\in\dom\Delta^{1/2}\cap M$ and note that
\[
\ip{\sum_ib_i\otimes y_i,\sum_jc_j\otimes z_j}_\Gamma
=\sum_{i,j}\tau(y_i^*\Gamma(b_i^*,c_j)z_j)
\]
is a positive semi-definite form on $(\dom\Delta^{1/2}\cap M)\otimes M$.
In particular, one has
\[
|\tau(x^*\Gamma(b^*,c)y)|\le\tau(x^*\Gamma(b^*,b)x)^{1/2}\tau(y^*\Gamma(c^*,c)y)^{1/2}.
\]
It follows that
\begin{align*}
\|\Gamma(b^*,c)\|_2
&=\sup\{|\tau(x^*\Gamma(b^*,c)y)| : x,y\in M,\ \|xx^*\|_2\le1,\|yy^*\|_2\le1\}\\
&\le\sup\{\tau(x^*\Gamma(b^*,b)x)^{1/2}\tau(y^*\Gamma(c^*,c)y)^{1/2} : \same\}\\
&\le\|\Gamma(b^*,b)\|_2^{1/2}\|\Gamma(c^*,c)\|_2^{1/2}\\
&\le4\|b\|_\infty^{1/2}\|\delta(b)\|^{1/2}\|c\|_\infty^{1/2}\|\delta(c)\|^{1/2}.
\end{align*}
\begin{lem}[Lemma 2.3 in \cite{peterson}]\label{lem:peterson}
For every $a,x\in M$, one has
\[
\|\zeta_\alpha(a)\tilde{\delta}_\alpha(x)
-\tilde{\delta}_\alpha(ax)\|
\le10\|x\|_\infty\|a\|_\infty^{1/2}\|\tilde{\delta}_\alpha(a)\|^{1/2}
\]
and
\[
\|\tilde{\delta}_\alpha(x)\zeta_\alpha(a)
-\tilde{\delta}_\alpha(xa)\|
\le10\|x\|_\infty\|a\|_\infty^{1/2}\|\tilde{\delta}_\alpha(a)\|^{1/2}.
\]
\end{lem}
\begin{proof}
One has
\begin{align*}
\zeta_\alpha(a)\tilde{\delta}_\alpha(x)
=\alpha^{-1/2}\delta(\zeta_\alpha(a)\zeta_\alpha(x))
-\tilde{\delta}_\alpha(a)\zeta_\alpha(x)
=:A_1-A_2.
\end{align*}
We note that $\|A_2\|\le\|x\|_\infty\|\tilde{\delta}_\alpha(a)\|$.
Let $\delta=V\Delta^{1/2}$ be the polar decomposition. Then, one has
\begin{alignat*}{3}
V^*A_1
&=\zeta_\alpha(a)\tilde{\Delta}_\alpha(x) &
&+\tilde{\Delta}_\alpha(a)\zeta_\alpha(x) &
&-\alpha^{-1/2}\Gamma(\zeta_\alpha(a),\zeta_\alpha(x))\\
&=:\qquad B_1 & &+\qquad B_2 & &-\qquad B_3
\end{alignat*}
in $L^2(M)$. We note that
$\|B_2\|\le\|x\|_\infty\|\tilde{\delta}_\alpha(a)\|$;
and by the estimate preceding to this lemma that
$\|B_3\|\le4\|x\|_\infty\|a\|_\infty^{1/2}\|\tilde{\delta}_\alpha(a)\|^{1/2}$.
Finally, one has
\begin{align*}
B_1
=\zeta_\alpha(a)\tilde{\Delta}_\alpha(x)
=\zeta_\alpha(a)(1-\theta_\alpha)(x)
\approx ax-\theta_\alpha(ax)
=\tilde{\Delta}_\alpha(ax).
\end{align*}
For the above estimates, we used
\[
\|\zeta_\alpha(a)x-ax\|_2\le\|x\|_\infty\|a-\zeta_\alpha(a)\|_2
\le\|x\|_\infty\|\tilde{\delta}_\alpha(a)\|_2
\]
and
\begin{align*}
\|\zeta_\alpha(a)\theta_\alpha(x)-\theta_\alpha(ax)\|_2
&\le\|x\|_\infty\|(\zeta_\alpha-\theta_\alpha)(a)\|_2
+\|\theta_\alpha(a)\theta_\alpha(x)-\theta_\alpha(ax)\|_2\\
&\le\|x\|_\infty(\|\tilde{\delta}_\alpha(a)\|_2+2\|a\|_\infty^{1/2}\|\tilde{\delta}_\alpha(a)\|^{1/2})
\end{align*}
(see Lemma \ref{lem:appcs}).
Consequently, one has
\[
\zeta_\alpha(a)\tilde{\delta}_\alpha(x)
\approx A_1\approx VB_1\approx\tilde{\delta}_\alpha(ax).
\]
This yields the first inequality.
Since the derivation is real, one obtains the second as well.
\end{proof}

We will need a vector-valued analogue of the above lemma.
Let
\[
\Omega=\{ \eta\in L^2(M\vt\bar{M}) :
(\id\otimes\bar{\tau})(\eta^*\eta)\le1\mbox{ and }
(\id\otimes\bar{\tau})(\eta\eta^*)\le1\}.
\]
We note that if $\{\xi_k\}$ is an orthonormal basis  of $L^2(M)$
and $\eta=\sum_{k=1}^\infty x_k\otimes\bar{\xi}_k$, then
$(\id\otimes\bar{\tau})(\eta^*\eta)=\sum_k x_k^*x_k$ and
$(\id\otimes\bar{\tau})(\eta\eta^*)=\sum_k x_kx_k^*$.
(These series converge \textit{a priori} in $L^1(M)$.)
We note that if $\eta\in\Omega$ and $b,c\in M$
with $\|b\|_\infty\|c\|_\infty\le1$,
then $\eta^*,(b\otimes1)\eta(c\otimes1)\in\Omega$.
\begin{lem}\label{lem:key}
For every $a\in M$ and $\eta\in\Omega$, one has
\[
\|(\zeta_\alpha(a)\otimes1)(\tilde{\delta}_\alpha\otimes1)(\eta)
-(\tilde{\delta}_\alpha\otimes1)((a\otimes1)\eta)\|_{\hh\vt L^2(\bar{M})}
\le20\|a\|_\infty^{1/2}\|\tilde{\delta}_\alpha(a)\|^{1/2}
\]
and
\[
\|(\tilde{\delta}_\alpha\otimes1)(\eta)(\zeta_\alpha(a)\otimes1)
-(\tilde{\delta}_\alpha\otimes1)(\eta(a\otimes1))\|_{\hh\vt L^2(\bar{M})}
\le20\|a\|_\infty^{1/2}\|\tilde{\delta}_\alpha(a)\|^{1/2}.
\]
\end{lem}

\begin{proof}
Let $a\in M$ be fixed and
define a linear map $T\colon M\to\hh$ by
\[
T(x)=\zeta_\alpha(a)\tilde{\delta}_\alpha(x)-\tilde{\delta}_\alpha(ax).
\]
By Lemma \ref{lem:peterson},
one has $\|T\|\le10\|a\|_\infty^{1/2}\|\tilde{\delta}_\alpha(a)\|^{1/2}$.
By the noncommutative little Grothendieck theorem
(Theorem 9.4 in \cite{pisier}),
there are states $f$ and $g$ on $M$ such that
\[
\|T(x)\|^2\le\|T\|^2\bigl(f(x^*x)+g(xx^*)\bigr)
\]
for all $x\in M$.
It follows that for $\eta=\sum_{k=1}^\infty x_k\otimes\bar{\xi}_k\in\Omega$,
one has
\begin{align*}
\|(\zeta_\alpha(a)\otimes1) & (\tilde{\delta}_\alpha\otimes1)(\eta)
-(\tilde{\delta}_\alpha\otimes1)((a\otimes1)\eta)\|_{\hh\vt L^2(\bar{M})}^2\\
&=\sum_k\|T(x_k)\|^2
\le\sum_k\|T\|^2\bigl(f(x_k^*x_k)+g(x_kx_k^*)\bigr)
\le2\|T\|^2.
\end{align*}
The second inequality follows similarly.
\end{proof}
\section{Proof of Theorems \ref{thm:A} and \ref{thm:B}}
Let $\G$ be a group and $(b,\pi,\hk)$ be a proper $1$-cocycle.
Replacing $(b,\pi,\hk)$ with
$(b\oplus\bar{b},\pi\oplus\bar{\pi},\hk\oplus\bar{\hk})$ and
considering an operator defined by $J_0(\xi\oplus\bar{\eta})=\eta\oplus\bar{\xi}$
if necessary, we may assume that there is a conjugate-linear involution $J_0$
on $\hk$ such that $J_0b(g)=b(g)$ and $J_0\pi_gJ_0=\pi_g$ for all $g\in\G$.
(Note that $\pi$ is amenable (resp.\ weakly sub-regular)
if and only if so is $\pi\oplus\bar{\pi}$.)

Let $M=Q\rtimes\G$ be the crossed product von Neumann algebra
of a finite von Neumann algebra $(Q,\tau)$ by a $\tau$-preserving
action $\sigma$ of $\G$.
We denote by $u_g$ the element in $M$ that corresponds to $g\in\G$.
We equip $\hh=L^2(Q)\otimes\ell^2(\G)\otimes\hk$ with
an $M$-$M$ bimodule structure by the following:
\[\begin{array}{rcccccc}
\hh & = & L^2(Q) & \otimes & \ell^2(\G)& \otimes & \hk\\
\mbox{left action by $g\in\G$} & : & \sigma_g & \otimes & \lambda_g & \otimes & \pi_g\\
\mbox{left action by $a\in Q$} & : & a & \otimes & 1 & \otimes & 1\\
\mbox{right action by $g\in\G$} & : & 1 & \otimes & \rho_g^{-1} & \otimes & 1\\
\mbox{right action by $a\in Q$} & : & {\displaystyle\sum_{h\in\G}\ \sigma_h(a)^{\op}} & \otimes & e_h & \otimes & 1
\end{array}\]
We define a conjugate-linear involution $J$ on $\hh$ by
\[
J(\widehat{a}\otimes\delta_g\otimes\xi)
=-\widehat{\sigma_{g^{-1}}(a^*)}\otimes\delta_{g^{-1}}\otimes J_0\pi_{g^{-1}}\xi,
\]
and the derivation $\delta\colon M\to\hh$ by
\[
\delta(au_g)=\widehat{a}\otimes\delta_g\otimes b(g)\in L^2(Q)\otimes\ell^2(\G)\otimes\hk
\]
for $a\in Q$ and $g\in\G$.
It is routine to check that $J$ intertwines the left and the right $M$-actions,
$J\delta(au_g)=\delta(\sigma_{g^{-1}}(a^*)u_{g^{-1}})=\delta((au_g)^*)$ and
moreover that $\delta$ is a real closable derivation satisfying
\[
\Delta(au_g)=\|b(g)\|^2au_g,\
\zeta_\alpha(au_g)=\sqrt{\frac{\alpha}{\alpha+\|b(g)\|^2}}au_g
\]
and
\[
\theta_\alpha(au_g)=\left(1-\sqrt{\frac{\|b(g)\|^2}{\alpha+\|b(g)\|^2}}\right)au_g
\]
for all $a\in Q$ and $g\in\G$.
In particular, all $\theta_\alpha$ belong to $C^*(Me_QM)$.

\begin{lem}\label{lem:coarse}
Suppose that $\pi$ is weakly contained in the regular representation.
Then, the $M$-$M$ bimodule $\hh$ is weakly contained in the
coarse bimodule $L^2(M)\vt L^2(M)$.
In particular, the left $M$-action on $\hh$ extends to
a u.c.p.\ map $\Psi\colon\B(L^2(M))\to\B(\hh)$
whose range commutes with the right $M$-action.
\end{lem}
\begin{proof}
It is well-known and not hard to see that if $\pi$ is weakly contained in the left
regular representation $\lambda$, then
the $M$-$M$ bimodule $\hh$ is weakly contained in
$\hat{\hh}:=L^2(Q)\vt\ell^2(\G)\vt\ell^2(\G)$, where
$(\pi,\hk)$ is replaced with $(\lambda,\ell^2(\G))$ in the definition of $\hh$.
Let $U$ be the unitary operator on $\hat{\hh}$
defined by
\[
U\widehat{a}\otimes\delta_h\otimes\delta_g
=\widehat{\sigma_g(a)}\otimes\delta_{gh}\otimes\delta_g.
\]
It is routine to check that
$U^*\lambda(M)U\subset\lambda(Q)\vt\C1\vt\B(\ell^2(\G))$ and
$U^*\rho(M^{\op})U\subset\rho(Q^{\op})\vt\B(\ell^2(\G))\vt\C1$,
where $\lambda$ and $\rho$ respectively stand for the left and right actions
on $\hat{\hh}$.
Since the ambient von Neumann algebras are amenable and commuting,
$\hat{\hh}$ and \textit{a fortiori} $\hh$
is weakly contained in the coarse $M$-$M$ bimodule, i.e.,
the binormal representation $\mu$ of $M\otimes M^{\op}$ on $\hh$ is
continuous w.r.t.\ the minimal tensor norm.
Hence, $\mu$ extends to a u.c.p.\ map $\tilde{\mu}$
from $\B(L^2(M))\vt M^{\op}$ into $\B(\hh)$.
We define $\Psi\colon\B(L^2(M))\to\B(\hh)$ by
$\Psi(x)=\tilde{\mu}(x\otimes1)$.
Since $M^{\op}$ is in the multiplicative domain of $\tilde{\mu}$,
the range of $\Psi$ commutes with the right $M$-action.
\end{proof}

For the following, let $P\subset M$ be an amenable von Neumann
subalgebra such that $P\not\preceq_M Q$, and $\GG\subset\Nor_M(P)$ be
a subgroup whose action on $P$ is weakly compact.
We may and will assume that $\U(P)\subset\GG$.
By definition, there exists a sequence $\eta_n\in L^2(M\vt\bar{M})_+$
such that
\begin{enumerate}
\item $\|\eta_n-(v\otimes\bar{v})\eta_n\|_2\to0$ for $v\in\U(P)$;
\item $\|\eta_n-\Ad(u\otimes\bar{u})(\eta_n)\|_2\to0$ for $u\in\GG$;
\item $\ip{(a\otimes1)\eta_n,\eta_n}=\tau(a)$ for all $a\in M$.
\end{enumerate}
We note that $\eta_n\in\Omega$.

\begin{lem}\label{lem:diff}
For every $\alpha>0$ and $a\in M$, one has
\[
\Lim_n\|(\tilde{\delta}_\alpha\otimes1)((a\otimes1)\eta_n)\|=\|a\|_2.
\]
\end{lem}
\begin{proof}
Note that $\|(a\otimes1)\eta_n\|_2=\|a\|_2$. Define a state on $\ip{M,e_Q}$ by
\[
\p_0(x)=\Lim_n\ip{(x\otimes1)\eta_n,\eta_n}.
\]
By construction,
$\p_0$ is a $P$-central state such that $\p_0|_M=\tau$.
Since $P\not\preceq_M Q$ and $\theta_\alpha a\in C^*(Me_QM)$,
Lemma \ref{lem:appcs} implies $\p_0(a^*\theta_\alpha^*\theta_\alpha a)=0$.
It follows that
\begin{align*}
\Lim_n\|(\tilde{\delta}_\alpha\otimes1)((a\otimes1)\eta_n)\|
&=\Lim_n\|((1-\theta_\alpha)a\otimes1)\eta_n\|_2\\
&=\Lim_n\|(a\otimes1)\eta_n\|_2\\
&=\|a\|_2.
\end{align*}
This completes the proof.
\end{proof}

For $\alpha>0$, a non-zero projection $p\in\GG'\cap M$ and $n$, we denote
\[
\eta_n^{p,\alpha}=(\tilde{\delta}_\alpha\otimes1)((p\otimes1)\eta_n)
\]
and define a state $\p_{p,\alpha}$ on $\B(\hh)\cap\rho(M^{\op})'$,
where $\rho(M^{\op})$ is the right $M$-action on $\hh$, by
\[
\p_{p,\alpha}(x)
=\|p\|_2^{-2}\Lim_n\ip{(x\otimes1)\eta_n^{p,\alpha},\eta_n^{p,\alpha}}.
\]

\begin{lem}\label{lem:appcent}
Let $a\in\GG''$. Then, one has
\[
\Lim_\alpha|\p_{p,\alpha}(\zeta_\alpha(a)x)-\p_{p,\alpha}(x\zeta_\alpha(a))|=0
\]
uniformly for $x\in\B(\hh)\cap\rho(M^{\op})'$ with $\|x\|_\infty\le1$.
\end{lem}
\begin{proof}
Let $u\in\GG$ and denote $u_\alpha=\zeta_\alpha(u)$.
By Lemma \ref{lem:key}, one has
\[
\Lim_n\|\eta_n^{p,\alpha}-(u_\alpha\otimes\bar{u})
 \eta_n^{p,\alpha}(u_\alpha\otimes\bar{u})^*\|
\le40\|\tilde{\delta}_\alpha(u)\|^{1/2}.
\]
Since $u_\alpha^*u_\alpha\le1$, one has for every $x\in(\rho(M^{\op})')_+$ that
\begin{align*}
\p_{p,\alpha}(u_\alpha^* x u_\alpha)
 &\geq\|p\|_2^{-2}\Lim_n\ip{(x\otimes1)(u_\alpha\otimes\bar{u})
 \eta_n^{p,\alpha}(u_\alpha\otimes\bar{u})^*,
 (u_\alpha\otimes\bar{u})\eta_n^{p,\alpha}(u_\alpha\otimes\bar{u})^*}\\
&\geq\p_{p,\alpha}(x)-80\|p\|_2^{-2}\|\tilde{\delta}_\alpha(u)\|^{1/2}\|x\|_\infty.
\end{align*}
By Lemma \ref{lem:stateapprox}, one obtains
\[
\|\p_{p,\alpha}(\,\cdot\,)-\p_{p,\alpha}(u_\alpha^*\,\cdot\,u_\alpha)\|
\le160\|p\|_2^{-2}\|\tilde{\delta}_\alpha(u)\|^{1/2}.
\]
In particular, $\Lim_\alpha\p_{p,\alpha}(1-u_\alpha^*u_\alpha)=0$ and
\[
\Lim_\alpha|\p_{p,\alpha}(u_\alpha x)-\p_{p,\alpha}(xu_\alpha)|=0
\]
uniformly for $x$ with $\|x\|_\infty\le1$.
This implies that
\[
\Lim_\alpha|\p_{p,\alpha}(\zeta_\alpha(a)x)-\p_{p,\alpha}(x\zeta_\alpha(a))|=0
\]
for each $a\in\lspan\GG$ and uniformly for $x\in\B(\hh)\cap\rho(M^{\op})'$
with $\|x\|_\infty\le1$.
However, by Lemma \ref{lem:key},
\begin{align*}
|\p_{p,\alpha}(x\zeta_\alpha(a))|
&=\|p\|_2^{-2}|\Lim_n\ip{(x\otimes1)(\zeta_\alpha(a)\otimes1)\eta_n^{p,\alpha},\eta_n^{p,\alpha}}|\\
&\le\|p\|_2^{-1}\|x\|_\infty\bigl(20\|a\|_\infty^{1/2}\|a\|_2^{1/2}+\|a\|_2),
\end{align*}
and likewise for $|\p_{p,\alpha}(\zeta_\alpha(a)x)|$.
Thus, by Kaplansky's Density Theorem, we are done.
\end{proof}

Now, we are in position to prove Theorems \ref{thm:A} and \ref{thm:B}.

\begin{proof}[Proof of Theorem \ref{thm:A}]
Let $\GG''=M$ and $\p_\alpha=\p_{1,\alpha}$.
By Lemma \ref{lem:appcent}, one has
\[
\Lim_\alpha|\p_\alpha(\zeta_\alpha(a)x)-\p_\alpha(x\zeta_\alpha(a))|=0
\]
for every $a\in \GG''=M$ and $x\in\B(\hh)\cap\rho(M^{\op})'$.
Since
\[
\|u_g-\zeta_\alpha(u_g)\|=1-\sqrt{\frac{\alpha}{\alpha+\|b(g)\|^2}}\to0
\]
as $\alpha\to\infty$, one has
\[
\Lim_\alpha|\p_\alpha(u_gxu_g^*)-\p_\alpha(x)|=0
\]
for every $g\in\G$ and $x\in\B(\hh)\cap\rho(M^{\op})'$.
Hence, the state $\p$ defined by
\[
\p(x)=\Lim_\alpha\p_\alpha(x)
\]
on $\B(\hk)\subset\B(\hh)\cap\rho(M^{\op})'$ is $\Ad\pi$-invariant.
Therefore, $\pi$ is an amenable representation, in contradiction to
the property {\HH}.
\end{proof}
\begin{proof}[Proof of Theorem \ref{thm:B}]
We use Haagerup's criterion for amenable von Neumann algebras (Lemma 2.2 in \cite{haagerup:dec}).
Let a non-zero projection $p\in\GG'\cap M$ and
a finite subset $F\subset\U(\GG'')$ be given arbitrary.
We need to show
\[
\|\sum_{u\in F} up\otimes \overline{up}\|_{M\vt\bar{M}}=|F|.
\]
Let $u\in\U(\GG'')$. By Lemmas \ref{lem:key} and \ref{lem:diff}, one has
\begin{align*}
\p_{p,\alpha}(\zeta_\alpha(up)^*\zeta_\alpha(up))
&=\|p\|_2^{-2}\Lim_n\|(\zeta_\alpha(up)\otimes1)
 (\tilde{\delta}_\alpha\otimes1)((p\otimes1)\eta_n)\|_2^2\\
&\geq\|p\|_2^{-2}\Lim_n\bigl(\|(\tilde{\delta}_\alpha\otimes1)((up\otimes1)\eta_n)\|_2
 -20\|\tilde{\delta}_\alpha(up)\|^{1/2}\bigr)^2\\
&\geq1-40\|p\|_2^{-1}\|\tilde{\delta}_\alpha(up)\|^{1/2}.
\end{align*}
Hence, by Lemma \ref{lem:appcent}, one has
\[
\Lim_\alpha|\p_{p,\alpha}(\zeta_\alpha(up)^*x\zeta_\alpha(up))-\p_{p,\alpha}(x)|=0
\]
uniformly for $x\in\B(\hh)\cap\rho(M^{\op})'$ with $\|x\|_\infty\le1$.
By Lemma \ref{lem:coarse}, the left $M$-action on $\hh$ extends to
a u.c.p.\ map $\Psi\colon\B(L^2(M))\to\B(\hh)\cap\rho(M^{\op})'$.
The state $\psi_{p,\alpha}=\p_{p,\alpha}\circ\Psi$ on $\B(L^2(M))$ satisfies
\[
\Lim_\alpha|\psi_{p,\alpha}(\zeta_\alpha(up)^*x\zeta_\alpha(up))-\psi_{p,\alpha}(x)|=0
\]
uniformly for $x\in\B(L^2(M))$ with $\|x\|_\infty\le1$.
By a standard convexity argument in cooperation with the Powers-St{\o}rmer inequality,
this implies that
\[
\Lim_\alpha\|\sum_{u\in F}\zeta_\alpha(up)\otimes\overline{\zeta_\alpha(up)}\|_{M\vt\bar{M}}
=|F|
\]
for the finite subset $F\subset\U(\GG'')$.
Since $\zeta_\alpha$ are u.c.p.\ maps, this yields
\[
\|\sum_{u\in F} up\otimes \overline{up}\|_{M\vt\bar{M}}\geq
\Lim_\alpha\|\sum_{u\in F}\zeta_\alpha(up)\otimes\overline{\zeta_\alpha(up)}\|_{M\vt\bar{M}}
=|F|.
\]
This completes the proof.
\end{proof}

\begin{proof}[Proof of Corollaries \ref{cor:A} and \ref{cor:B}]
The corollaries follow from the corresponding Theorems
and Lemma \ref{conjugatecartan}, because all the von Neumann algebras
in consideration have the complete
metric approximation property and hence all amenable subalgebras are
weakly compact (Theorem \ref{thm:cmapwc}).
\end{proof}

\begin{rem}
With the same argument as above one can show the following. Let
$(M,\tau)$ be a finite von Neumann algebra having a real closable
derivation $(\delta,\hh)$ such that (1) $\delta^*\delta$ has compact
resolvent; and (2) $\hh$ is weakly contained in the coarse bimodule.
Then, for every diffuse amenable von Neumann subalgebra $P\subset M$
and any subgroup $\GG\subset\Nor_M(P)$ whose action on $P$ is weakly
compact, the von Neumann subalgebra $\GG''$ is amenable.
\end{rem}
\section{Cocycle rigidity}
We fix a notation for profinite actions.
An action $\G\curvearrowright (X,\mu)$ is said to be \emph{profinite}
if $(X,\mu)$ is the projective limit of finite-cardinality probability
spaces $(X_n,\mu_n)$ on which $\G$ acts consistently.
We will identify $L^\infty(X_n,\mu_n)$ as the corresponding
$\G$-invariant finite-dimensional von Neumann subalgebra of
$L^\infty(X,\mu)$. The same thing for $L^2$.
We write $X=\bigsqcup_a X_{a,n}$ for the partition of $X$ corresponding
to $X_n$, i.e., the characteristic functions of $X_{a,n}$'s are
the non-zero minimal projections in $L^\infty(X_n)$.

\begin{defn}
Let $\pi\colon\G\curvearrowright\hh$ be a unitary representation.
We say $\pi$ \emph{has a spectral gap}
if there are a finite subset $F\subset\G$ and $\kappa>0$, called
\emph{a critical pair}, satisfying the following property:
denoting by $P$ the orthogonal projection
of $\hh$ onto the subspace of $\pi$-invariant vectors, one has
\[
\kappa\|\xi-P\xi\|\le\max_{g\in F}\|\xi-\pi_g\xi\|
\]
for every $\xi\in\hh$. (This is equivalent to that the point $1$ is
isolated (if it exists) in the spectrum of the self-adjoint operator
$(2|F|)^{-1}\sum_{g\in F}(\pi_g+\pi_g^*)$ on $\hh$.) We say that
$\pi$ \emph{has a stable spectral gap} if the unitary representation
$\pi\otimes\bar{\pi}$ of $\G$ on $\hh\vt\bar{\hh}$ has a spectral
gap. (Note that we allow $\rank P\geq1$.)

When the unitary representation $\pi$ arises from a p.m.p.\ action
$\G\curvearrowright X$, we simply say \emph{$\G\curvearrowright X$
has a (stable) spectral gap} if $\pi$ has. Assume moreover that the
action $\G\curvearrowright X$ is profinite. We say
\emph{$\G\curvearrowright X$ has a stable spectral gap with growth
condition} if there are a critical pair $(F,\kappa)$ such that
$\pi^F$, the restriction of $\pi$ to the subgroup of $\G$ generated
by $F$, does not have a subrepresentation of infinite multiplicity.
\end{defn}

Suppose that $\G\curvearrowright \varprojlim X_n$ has a stable spectral gap.
Then, $\pi^F$ has finitely many equivalence classes of irreducible
subrepresentations of any given dimension $k\in\N$. (See \cite{hrv}.)
It follows that the growth condition is equivalent to that the minimal
dimension $k_n$ of a non-zero subrepresentation
of $\pi^F|_{L^2(X)\ominus L^2(X_n)}$ tends to infinity.
\begin{lem}\label{lem:ssgg}
Let $\G\curvearrowright X$ be a p.m.p.\ action
which is profinite and has a stable spectral gap with growth condition.
Let $F\subset\G$ and $\kappa>0$ be a critical pair.
Then, for any $k\in\N$ and unitary elements
$\{u_g\}_{g\in F}$ on the $k$-dimensional Hilbert space $\ell^2_k$,
one has
\[
\frac{\kappa^2}{2}(1-\frac{k}{k_n})\|\xi-P_{L^2(X_n)\vt\ell^2_k}\xi\|_2
\le\max_{g\in F}\|\xi-(\pi_g\otimes u_g)\xi\|_2
\]
for every $\xi\in L^2(X)\vt\ell^2_k$ and $n\in\N$.
\end{lem}
\begin{proof}
We denote $L^2(X_n)^\perp=L^2(X)\ominus L^2(X_n)$.
It suffices to show
\[
\frac{\kappa^2}{2}(1-\frac{k}{k_n})\|\xi\|_2
\le\max_{g\in F}\|\xi-(\pi_g\otimes u_g)\xi\|_2
\]
for $\xi\in L^2(X_n)^\perp\vt\ell^2_k$.
We assume $\|\xi\|_2=1$ and denote the right hand side of the asserted
inequality by $\e$.
We view $\xi$ as a Hilbert-Schmidt
operator $T_\xi$ from $\overline{\ell^2_k}$ into $L^2(X_n)^\perp$.
Note that
\[
\|T_\xi-\pi_g T_\xi\bar{u}_g^*\|_{2}=\|\xi-(\pi_g\otimes u_g)\xi\|_2\le\e.
\]
Hence by the Powers-St{\o}rmer inequality, the Hilbert-Schmidt operator
$S_\xi=(T_\xi T_\xi^*)^{1/2}$ on $L^2(X_n)^\perp$ satisfies
\begin{align*}
\| S_\xi-\pi_g S_\xi \pi_g^*\|_2^2
&\le\| T_\xi T_\xi^*-\pi_g T_\xi T_\xi^* \pi_g^*\|_1\\
&\le\|T_\xi+\pi_g T_\xi\bar{u}_g^*\|_{2}\|T_\xi-\pi_g T_\xi\bar{u}_g^*\|_{2}\\
&\le2\e.
\end{align*}
By the stable spectral gap property, one has
\[
\| S_\xi- P(S_\xi)\|_2^2\le2\e/\kappa^2.
\]
Since $P(S_\xi)$ commutes with $\pi_g$ for all $g\in F$,
growth condition implies that
$P(S_\xi)=\sum_i\gamma_ir_i^{-1/2}Q_i$ for some $\gamma_i\geq0$ and
mutually orthogonal projections $Q_i$ with $r_i=\Tr(Q_i)\geq k_n$.
Since $S_\xi$ has rank at most $k$,
denoting its range projection by $R$, one has
\begin{align*}
\|P(S_\xi)\|_2^2 &=\ip{S_\xi, P(S_\xi)}\\
&=\sum_i\gamma_ir_i^{-1/2}\Tr(Q_iRS_\xi Q_i)\\
&\le\bigl(\sum_i\gamma_i^2r_i^{-1}\Tr(Q_iRQ_i)\bigr)^{1/2}
 \bigl(\sum_i\Tr(Q_iS_\xi^*S_\xi Q_i)\bigr)^{1/2}\\
&\le\bigl(\sum_i\gamma_i^2k_n^{-1}k\bigr)^{1/2}\|S_\xi\|_2\\
&=(k/k_n)^{1/2}\|P(S_\xi)\|_2.
\end{align*}
By combining two inequalities, one obtains
\[
1-(k/k_n)\le1-\|P(S_\xi)\|_2^2=\| S_\xi- P(S_\xi)\|_2^2\le2\e/\kappa^2
\]
and hence the desired inequality.
\end{proof}

Recall that a cocycle of $\G\curvearrowright X$ with values in a
group $\Lambda$ is a measurable map $\alpha\colon\G\times X\to\Lambda$
satisfying the cocycle identity:
\[
\forall g,h\in\G,\,\mbox{$\mu$-a.e. }x\in X,\quad
\alpha(g,hx)\alpha(h,x)=\alpha(gh,x).
\]
A cocycle $\alpha$ which is independent of the $x$-variable is
said to be \emph{homomorphism} for the obvious reason.
Cocycles $\alpha$ and $\beta$ are said to be \emph{equivalent} if there is a measurable
map $\phi\colon X\to\Lambda$ such that $\beta(g,x)=\phi(gx)\alpha(g,x)\phi(x)^{-1}$
for all $g\in\G$ and $\mu$-a.e.\ $x\in X$.
\begin{lem}\label{lem:cocfg}
Let $\G=\G_1\times\G_2$ and $\G\curvearrowright X=\varprojlim X_n$ be a
p.m.p.\ profinite action such that
$\G_2\curvearrowright X$ has a stable spectral gap with growth condition.
Let $(N,\tau)$ be a finite type $\mathrm{I}$ von Neumann algebra, and
$\alpha\colon \G\times X\to\U(N)$ be a cocycle.
Then, for every $\e>0$, there exists $n\in\N$
such that
\[
\int_X \|\alpha(g,x)-\alpha'_{a(x),n}(g)\|_2^2\,dx\le\e
\]
for all $g\in\ker(\G_1\to\Aut(X_n))$,
where $\alpha'_{a,n}(g)=|X_{a,n}|^{-1}\int_{X_{a,n}}\alpha(g,y)\,dy$
and $a(x)$ is such that $x\in X_{a(x),n}$.
\end{lem}
\begin{proof}
It suffices to consider each direct summand of $N$ and
hence we may assume $N=\M_k(\C)\otimes A$, where
$A$ is an abelian von Neumann algebra.
For every $g\in\G$, we define $w_g\in L^\infty(X)\vt N=L^\infty(X,N)$ by
$w_g(x)=\alpha(g,g^{-1}x)$.
Then, it becomes a unitary $1$-cocycle for
$\tilde{\sigma}=\sigma\otimes\id_N$:
\[
w_{gh}=w_g\tilde{\sigma}_g(w_h).
\]
Let $F\subset\G_2$ and $\kappa>0$ be a critical pair
for the stable spectral gap of $\G_2\curvearrowright X$.
Let $\delta=\e\kappa^2/8$ and
take $m\in\N$ and unitary elements $w_h'\in L^\infty(X_m)\vt N$
such that $\| w_h-w_h' \|_2<\delta$ for every $h\in F$.
For the rest of the proof, we fix $g\in\ker(\G_1\to\Aut(X_m))$.
Since $w_h'=\tilde{\sigma}_g(w_h')$, one has
\[
w_gw_h'\approx w_g\tilde{\sigma}_g(w_h)=w_{gh}
=w_{hg}\approx w_h'\tilde{\sigma}_h(w_g),
\]
for every $h\in F$.
We define trace-preserving $*$-automorphisms $\pi_h$
on $L^\infty(X)\vt N$ by
\[
\pi_h(x)=\Ad(w_h')\circ\tilde{\sigma}_h(x)
\]
and note that $\|w_g-\pi_h(w_g)\|_2\le2\delta$ for every $h\in F$.
We write $\tilde{\pi}_h$ for the restriction of
$\pi_h$ to $L^\infty(X_m)\vt N$.
Note that $\tilde{\pi}_h$ acts as identity on $\C1\vt A\subset L^\infty(X_m)\vt N$.

Let $\{p_a\}$ be the set of non-zero minimal projections in $L^\infty(X_m)$
and define an isometry $V\colon L^2(X)\to L^2(X)\vt L^2(X_m)$ by
$V\xi=|X_m|^{1/2}\sum_a p_a\xi\otimes p_a$. (Here $|X_m|$
stands for the cardinality of the atoms of $X_m$.)
We claim that $(V\otimes1)\pi_h=(\sigma_h\otimes\tilde{\pi}_h)(V\otimes1)$.
Indeed, if $w_h'=\sum_a \sigma_h(p_a)\otimes y_a$, then
\begin{align*}
(\sigma_h\otimes\tilde{\pi}_h)(V\otimes1)(\xi\otimes c)
&=|X_m|^{1/2}(\sigma_h\otimes\tilde{\pi}_h)\sum_a p_a\xi\otimes p_a\otimes c\\
&=|X_m|^{1/2}\sum_a\sigma_h(p_a\xi)\otimes\sigma_h(p_a)\otimes y_a c y_a^*\\
&=V\sum_a\sigma_h(p_a\xi)\otimes y_a c y_a^*\\
&=V\pi_h(\xi\otimes c)
\end{align*}
for all $\xi\in L^2(X)$ and $c\in L^2(N)$.
Now, it follows that
\[
\max_{h\in F}\|(V\otimes1)w_g-(\sigma_h\otimes\tilde{\pi}_h)(V\otimes1)w_g\|_2\le2\delta.
\]
We observe that if $\tilde{\pi}_h$ is viewed as a unitary operator on
$L^2(X_m)\vt L^2(\M_k(\C))\vt L^2(A)$, then it lives in $\B(L^2(X_m)\vt L^2(\M_k(\C)))\vt A$.
Hence Lemma \ref{lem:ssgg} applies and one obtains
\[
\frac{\kappa^2}{2}(1-\frac{mk^2}{k_n})
 \|(V\otimes1)w_g-(P_{L^2(X_n)}\otimes1\otimes1)(V\otimes1)w_g\|_2
\le2\delta
\]
for every $n\in\N$.
Finally take $n$ to be such that $n\geq m$ and $k_n\geq2mk^2$.
Since $(P_{L^2(X_n)}\otimes1)V=VP_{L^2(X_n)}$ for $n\geq m$,
one has
\begin{align*}
(\int_X \|\alpha(g,x)-\alpha'_{a(x),n}(g)\|_2^2\,dx)^{1/2}
&=\|w_g-(P_{L^2(X_n)}\otimes1)w_g\|_2\\
&\le4\delta/(\kappa^2(1-(\frac{mk^2}{k_n})^{1/2}))\le\e.
\end{align*}
We note that $\ker(\G_1\to\Aut(X_n))\subset\ker(\G_1\to\Aut(X_m))$.
\end{proof}

We combine the above result with results of Ioana in
\cite{ioana:cocycle}, to obtain the following cocycle rigidity
result for profinite actions of product groups.

\begin{thm}\label{thm:cocyuntf}
Let $\G=\G_1\times\G_2$ and $\G\curvearrowright X$ be an ergodic
p.m.p.\ profinite action such that $\G_i\curvearrowright X$ has a
stable spectral gap with growth condition, for each $i=1,2$. Let
$\Lambda$ be a finite group and $\alpha\colon\G\times X\to\Lambda$
be a cocycle. Then, there exists a finite index subgroup
$\G'\subset\G$ such that for each $\G'$-ergodic component $X'\subset
X$, the restricted cocycle $\alpha|_{\G'\times X'}$ is equivalent to
a homomorphism from $\G'$ into $\Lambda$.
\end{thm}
\begin{proof}
The proof of this theorem is very similar
to that of Theorem B in \cite{ioana:cocycle},
and hence it will be rather sketchy.
Let $Z=X\times X\times\Lambda$ and we will consider
the unitary representation $\pi\colon\G\curvearrowright L^2(Z)$
induced by the m.p.\ transformation
\[
g(x,y,t)=(gx,gy,\alpha(g,x)t\alpha(g,y)^{-1}).
\]
Let $\e>0$ be arbitrary.
Since $\Lambda$ is discrete, Lemma \ref{lem:cocfg} implies that
there are a normal finite index subgroup $\G'$ and $n\in\N$ such that
$\|\pi(g)\xi_n-\xi_n\|_2<\e$ for all $g\in\G'$, where
$\xi_n=|X_n|^{1/2}\sum_a\chi_{X_{a,n}\times X_{a,n}\times\{e\}}$.
It follows that the circumcenter of $\pi(\G')\xi_n$ is
a $\pi(\G')$-invariant vector which is close to $\xi_n$.
Since $\G\curvearrowright X$ is ergodic and
$\G'$ is a normal finite index subgroup in $\G$,
there are a $\G'$-ergodic component $X'\subset X$ and
a finite subset $E\subset\G$ such that $X=\bigsqcup_{s\in E}sX'$.
Thus, there are $\G'$-ergodic components $X'_1,X'_2\subset X$ such that
$\xi'=|X'|^{-1}\chi_{X'_1\times X'_2\times\{e\}}$ is
close to a $\pi(\G')$-invariant vector.
We may assume that $X'_1=X'$.
By Corollary 2.2 in \cite{ioana:cocycle}, the cocycle $\alpha|_{\G'\times X'}$
is equivalent to a homomorphism $\theta$ via $\phi\colon X'\to\Lambda$, i.e.,
$\theta(g)=\phi(gx)\alpha(g,x)\phi(x)^{-1}$.
We observe that $\alpha|_{\G'\times sX'}$ is equivalent to $\theta\circ\Ad(s^{-1})$.
Indeed, one has
\begin{align*}
\theta(s^{-1}gs)
&=\phi(s^{-1}gsx)\alpha(s^{-1}gs,x)\phi(x)^{-1}\\
&=\phi(s^{-1}gsx)\alpha(s^{-1},gsx)\alpha(g,sx)\alpha(s,x)\phi(x)^{-1}\\
&=\psi(gsx)\alpha(g,sx)\psi(sx)^{-1},
\end{align*}
where $\psi(sx)=\phi(x)\alpha(s^{-1},sx)$ for $s\in E$ and $x\in X'$.
\end{proof}

\begin{proof}[Proof of Theorem \ref{thm:D}]
By Theorem B and Remark 3.1 in \cite{ioana:cocycle},
it suffices to show that the unitary representation
$\pi\colon\G\curvearrowright L^2(X\times X\times\Lambda)$ has a spectral gap.
Let $\Lambda_j$ be the finite quotients of $\Lambda$.
Since $\Lambda$ is residually finite the unitary representation
$\pi$ is weakly contained in the direct sum $\bigoplus\pi_j$, where
$\pi_j$ is the unitary representation induced by
$\G\curvearrowright X\times X\times\Lambda_j$ using the same $1$-cocycle
composed with the quotient $\Lambda\to\Lambda_j$.
Thus, it suffices to show that $\pi_j$'s have a uniform spectral gap.
We prove this by showing that each $\pi_j$ is contained in
a direct sum of finite representations; then the uniformity follows
from property $(\tau)$.
We may assume that $\Lambda$ is finite and $\pi_j=\pi$.
By Theorem \ref{thm:cocyuntf}, there is a finite index subgroup $\G'$
such that for each $\G'$-ergodic component $X'_k\subset X$,
the restricted cocycle $\alpha|_{\G'\times X'_k}$ is equivalent
to a homomorphism $\theta_k\colon\G'\to\Lambda$
via $\phi_k\colon X'_k\to\Lambda$, i.e.,
$\theta_k(g)=\phi_k(gx)\alpha(g,x)\phi_k(x)^{-1}$.
Let $\sigma_{k,l}$ be the automorphism on $X'_k\times X'_l\times\Lambda$
defined by $\sigma_{k,l}(x,y,t)=(x,y,\phi_k(x)t\phi_l(y)^{-1})$.
Then, for $g\in\G'$, one has
\begin{align*}
\sigma_{k,l}g\sigma_{k,l}^{-1}(x,y,t)
&=(gx,gy,\phi_k(gx)\alpha(g,x)\phi_k(x)^{-1}t\phi_l(y)\alpha(g,y)^{-1}\phi_l(gy)^{-1})\\
&=(gx,gy,\theta_k(g)t\theta_l(g)^{-1}).
\end{align*}
Since $X'$ is profinite, this implies that
the unitary representation $\pi|_{\G'}$ is contained in a direct sum of
finite representations (of the form $\G'\curvearrowright X_n\times X_n\times \Lambda$,
$g(x,y,t)=(gx,gy,\theta_k(g)t\theta_l(g)^{-1})$).
Since $\G'$ has finite index in $\G$, the unitary representation
$\pi\subset\Ind_{\G'}^{\G}(\pi|_{\G'})$
is contained in a direct sum of
finite representations. This completes the proof.
\end{proof}

The following two lemmas are well-known, but we include the proof
for the reader's convenience.

\begin{lem}
Let $\G\geq\Delta_1\geq\Delta_2\geq\cdots$ be a decreasing sequence of
finite index normal subgroups.
Then, the left-and-right action
$\G\times\G\curvearrowright \varprojlim\G/\Delta_n$
is essentially-free if and only if $\lim_n|Z_n(g)|/|\G/\Delta_n|=0$
for every $g\in\G$ with $g\neq e$,
where $Z_n(g)$ is the centralizer group of $g$ in $\G/\Delta_n$.
\end{lem}
\begin{proof}
The `only if' part is trivial. We prove the `if' part.
Note that the condition implies that $\bigcap\Delta_n=\{e\}$.
Let $(g,h)\in\G\times\G$ and observe that
\[
|\{ x\in\varprojlim\G/\Delta_n : (g,h)x=x\}|
=\lim_n\frac{|\{x\in\G/\Delta_n : gxh^{-1}=x\}|}{|\G/\Delta_n|}.
\]
If $g=e$, then $(g,h)$ acts freely unless $h=e$, too.
Thus, let $g\neq e$. If $x,y\in\G/\Delta_n$ are such that
$gxh^{-1}=x$ and $gyh^{-1}=y$, then one has
$gxy^{-1}g^{-1}=xy^{-1}$, i.e., $xy^{-1}\in Z_n(g)$.
It follows that $|\{x\in\G/\Delta_n : gxh^{-1}=x\}|\le|Z_n(g)|$.
\end{proof}

\begin{lem}\label{lem:pslcent}
Let $F$ be a finite field.
Then, for every $g\in\PSL(2,F)$ with $g\neq e$,
one has $|Z(g)|/|\PSL(2,F)|\le2/(|F|-1)$.
\end{lem}
\begin{proof}
Since the characteristic polynomial of $g$ is quadratic,
it can be factorized in some quadratic extension $\tilde{F}$ of $F$.
Thus $g$ is conjugate to a Jordan normal form in $\PSL(2,\tilde{F})$.
Now, it is not hard to see that the centralizer of $g$ in $\PSL(2,\tilde{F})$
has cardinality at most $|\tilde{F}|=|F|^2$.
On the other hand, it is well-known that $|\SL(2,F)|=|F|(|F|^2-1)$.
\end{proof}

\begin{proof}[Proof of Corollary \ref{cor:D}]
Since $\SL(2,\Z[\sqrt{2}])$ is an irreducible lattice in
$\SL(2,\R)^2$, it has property $(\tau)$ (see Section 4.3 in
\cite{lubotzky}) and the property {\HH}$^+$ (cf. Theorem 2.3). By
the above lemmas, the action $\G\curvearrowright X$ is
essentially-free. Indeed, consider the homomorphism from
$\PSL(2,(\Z/p_1\cdots p_n\Z)[\sqrt{2}])$ onto $\PSL(2,F)$, where $F$
is the field either $\Z/p_n\Z$ or $(\Z/p_n\Z)[\sqrt{2}]$, depending
on whether the equation $x^2=2$ is solvable in $\Z/p_n\Z$ or not;
and apply Lemma \ref{lem:pslcent} at $\PSL(2,F)$. Therefore, by
Corollary \ref{cor:A}, $L^\infty(X)$ is the unique Cartan subalgebra
of $L^\infty(X)\rtimes\G$. It follows that the isomorphism of von
Neumann algebras $L^\infty(X)\rtimes\G\cong
(L^\infty(Y)\rtimes\Lambda)^t$ gives rise to a stable orbit
equivalence between $\G\curvearrowright X$ and
$\Lambda\curvearrowright Y$. The growth condition of Theorem
\ref{thm:D} is satisfied because $p_k$'s are mutually distinct
primes and $\PSL(2,(\Z/p_1\cdots
p_n\Z)[\sqrt{2}])\cong\prod\PSL(2,(\Z/p_k\Z)[\sqrt{2}])$. Therefore,
Theorem \ref{thm:D} is applicable to the orbit equivalence cocycle
$\alpha\colon \G\times X\to\Lambda$. For the rest of the proof, see
\cite{ioana:cocycle}.
\end{proof}
\section{$\mathrm{II}_1$-factors with more than one Cartan subalgebra}\label{sec:more}
\begin{proof}[Proof of Theorem \ref{thm:C}]
Since $H\curvearrowright^\sigma X$ is an ergodic and profinite action,
one has $X=\varprojlim H/H_n$ for some decreasing sequence
$H=H_0\supset H_1\supset\cdots$ of finite-index subgroups of $H$
such that $\bigcap H_n=\{e\}$.
Recall that a function $f$ is called an eigenfunction of $H$ if
there is a character $\chi$ on $H$ such that $\sigma_h(f)=\chi(h)f$
for every $h\in H$. We observe that every unitary eigenfunction normalizes
$L(H)$ in $L^\infty(X)\rtimes H$, and that $L^\infty(X)$
is spanned by unitary eigenfunctions since $L^\infty(H/H_n)$
is spanned by characters.
This proves that $L(H)$ is regular in $L^\infty(X)\rtimes\G$.
To prove that $L(H)$ is maximal abelian, let
$a\in L(H)'\cap L^\infty(X)\rtimes\G$ be given and
$a=\sum_{g\in\G} a_gu_g$ be the Fourier expansion.
Then, $[a,u_h]=0$ implies $\sigma_h(a_g)=a_{hgh^{-1}}$ for
all $g\in\G$ and $h\in H$.
In particular, one has $\|a_{hgh^{-1}}\|_2=\|a_g\|_2$.
Since $\sum_g\|a_g\|_2^2=\|a\|_2^2<\infty$,
the relative ICC condition implies that $a_g=0$ for all $g\not\in H$.
But for $g\in H$, ergodicity of $H\curvearrowright X$ implies
that $a_g\in\C1$. This proves $a\in L(H)$.

For the second assertion, recall that weak compactness is
an invariant of a Cartan subalgebra (Proposition 3.4 in \cite{I}).
We prove that $L(H)$ is not weakly compact in $L^\infty(X)\rtimes\G$.
Suppose by contradiction that it is weakly compact.
Then, by Proposition 3.2 in \cite{I}, there is a state $\p$ on
$\B(\ell^2(H))$ which is invariant under the $H\rtimes\G$-action.
Restricting $\p$ to $\ell^\infty(H)$, we obtain an $H\rtimes\G$-invariant mean.
This contradicts the assumption.
\end{proof}
Corollary \ref{cor:c} is an immediate consequence of Theorem \ref{thm:C}.
Here we give another example for which Theorem \ref{thm:C} applies.
Let $K$ be a residually-finite additive group such that $|K|>1$, and
$\G_0$ be a residually-finite non-amenable group.
The wreath product $\G=K\wr\G_0$ is defined to be the semidirect product
of $H=\bigoplus_{\G_0} K$ by the shift action of $\G_0$.
Then, there is a decreasing sequence $H_0\supset H_1\supset\cdots$
of $\G_0$-invariant finite-index subgroups of $H$ such that $\bigcap H_n=\{0\}$.
Indeed, let $K_0\supset K_1\supset\cdots$
(resp.\ $\G_{0,0}\supset\G_{0,1}\supset\cdots$)
be finite-index subgroups of $K$ (resp.\ $\G_0$)
such that $\bigcap K_n=\{0\}$ (resp.\ $\bigcap\G_{0,n}=\{e\}$).
Then, the ``augmentation subgroups''
\[
H_n =\{ (a_g)_{g\in\G_0}\in H
 : \sum_{h\in\G_{0,n}}a_{gh}\in K_n\mbox{ for all $g\in\G_0$}\},
\]
which is the kernel of the homomorphism onto
$\bigoplus_{\G_0/\G_{0,n}}K/K_n$, satisfy the required conditions.
It follows from Theorem \ref{thm:C} that the $\mathrm{II}_1$-factor
\[
L^\infty(\varprojlim H/H_n)\rtimes\G
\]
has two non-conjugate Cartan subalgebras, namely $L(H)$ and
$L^\infty(\varprojlim H/H_n)$.

\end{document}